\documentclass{amsart}
\usepackage{amsfonts}
\usepackage{amsmath}
\usepackage{amssymb}
\usepackage{amsthm}
\newtheorem{theorem}{Theorem}
\newtheorem{prop}[theorem]{Proposition}
\newtheorem{lemma}[theorem]{Lemma}

\newtheorem{exmp}[theorem]{Example}
\newtheorem{cor}[theorem]{Corollary}
\newtheorem{fact}[theorem]{Fact}

\begin{document}

\author{Mark Pankov}
\title{Base subsets of Symplectic Grassmannians}

\address{Department of Mathematics and Information Technology,
University of Warmia and Mazury, {\. Z}olnierska 14A, 10-561 Olsztyn, Poland}
\email{pankov@matman.uwm.edu.pl}

\begin{abstract}
Let $V$ and $V'$ be $2n$-dimensional vector spaces over fields $F$ and $F'$.
Let also $\Omega: V\times V\to F$ and $\Omega': V'\times V'\to F'$
be non-degenerate symplectic forms.
Denote by $\Pi$ and $\Pi'$ the associated $(2n-1)$-dimensional projective spaces.
The sets of $k$-dimensional totally isotropic subspaces of $\Pi$ and $\Pi'$
will denoted by ${\mathcal G}_{k}$ and ${\mathcal G}'_{k}$, respectively.
Apartments of the associated buildings intersect ${\mathcal G}_{k}$ and ${\mathcal G}'_{k}$
by so-called base subsets.
We show that every mapping of ${\mathcal G}_{k}$ to ${\mathcal G}'_{k}$
sending base subsets to base subsets is induced by a symplectic embedding of
$\Pi$ to $\Pi'$.
\end{abstract}

\maketitle

\section{Introduction}
An incidence geometry of the rank $n$ has the following ingredients:
a set ${\mathcal G}$ whose elements are called {\it subspaces},
a symmetric {\it incidence relation} on ${\mathcal G}$,
and a surjective {\it dimension function}
$$\dim: {\mathcal G}\to \{0,1,\dots,n-1\}$$
such that the restriction of
this function to every maximal flag is bijective
(flags are set of mutually incident subspaces).

A {\it Tits building} \cite{Tits} is an incidence geometry together with a
family of isomorphic subgeometries called {\it apartments} and
satisfying a certain collection of axioms.
One of these axioms
says that for any two flags there is an apartment containing them.

Let us consider an incidence geometry of the rank $n$
whose set of subspaces is denoted by ${\mathcal G}$.
For every $k\in \{0,1,\dots, n-1\}$ we
denote by ${\mathcal G}_{k}$ the {\it Grassmannian} consisting of
all $k$-dimensional subspaces.
If this geometry is a building then
the intersection of ${\mathcal G}_{k}$ with an apartment
is called {\it the shadow} of this apartment in ${\mathcal G}_{k}$
\cite{Cohen}.
In the projective and symplectic cases the intersections
of apartments with Grassmannians
are known as {\it base subsets} \cite{Pankov1,Pankov2,Pankov3}.

Let $f$ be a bijective transformation of ${\mathcal G}_{k}$
preserving the family of the shadows of apartments.
It is natural to ask: can $f$ be extended to an automorphism of the corresponding geometry?
This problem  was solved in \cite{Pankov1} for buildings of the type $A_{n}$,
in this case $f$ is induced by a collineation of the associated projective space
to itself or the dual projective space
(the second possibility can be realized only for the case when $n=2k+1$).
A more general result can be found in \cite{Pankov2}.

In the present paper we show that the extension is possible
for symplectic buildings.

Note that apartment preserving transformations of the chamber set
(the set of maximal flags) of a spherical building
are induced by automorphisms of the corresponding complex;
this follows from the results given in \cite{AVM}.

\section{Symplectic geometry}
Let $V$ be a $2n$-dimensional vector space over a field $F$, and
let $$\Omega: V\times V\to F$$ be a non-degenerate symplectic form.
Denote by $\Pi=(P,{\mathcal L})$
the $(2n-1)$-dimensional projective space associated with $V$
(points are $1$-dimensional subspaces of $V$ and lines are defined by
$2$-dimensional subspaces).

We say that two points $p,q\in P$ are {\it orthogonal} and write $p\perp q$ if
$$p=<x>,\;q=<y>\;\mbox{ and }\;\Omega(x,y)=0.$$
Similarly, two subspaces $S$ and $U$ of $\Pi$ will be called {\it orthogonal} ($S\perp U$) if
$p\perp q$ for any $p\in S$ and $q\in U$.
The orthogonal complement to a
subspace $S$ (the maximal subspace orthogonal to $S$) will be
denoted by $S^{\perp}$, if $S$ is $k$-dimensional then the
dimension of $S^{\perp}$ is equal to $2n-k-2$.

A base $\{p_{1},\dots,p_{2n}\}$ of $\Pi$ is said to be {\it symplectic}
if for each $i\in \{1,\dots,2n\}$ there exists
unique $\sigma(i)\in \{1,\dots,2n\}$ such that
$$p_{i}\not\perp p_{\sigma(i)}.$$
($p_{i}$ and $p_{\sigma(i)}$ are non-orthogonal).

A subspace $S$ of $\Pi$ is called {\it totally isotropic} if any
two points of $S$ are orthogonal; in other words, $S\subset S^{\perp}$.
The latter inclusion implies that the dimension of a
totally isotropic subspace is not greater than $n-1$.

Now consider the incidence geometry of totally isotropic subspaces.
For every symplectic base $B$ the subgeometry consisting of all
totally isotropic subspaces spanned by points of $B$ is
{\it the symplectic apartment associated with} $B$.
It is well-known that the incidence geometry of totally isotropic subspaces
together with the family of all symplectic apartments is a building of the type $C_{n}$.

We write ${\mathcal G}_{k}$ for the set of all $k$-dimensional  totally isotropic subspaces.
The set of all $k$-dimensional totally isotropic
subspaces spanned by points of a symplectic base will be called
the {\it base subset} of ${\mathcal G}_{k}$
associated with (defined by) this base.

\begin{prop}\label{pr1}
Every base subset of ${\mathcal G}_{k}$ consists of
$$2^{k+1}\binom{n}{k+1}$$
elements.
\end{prop}

\begin{proof}
Let $B=\{p_{1},\dots, p_{2n}\}$ be a symplectic base and
${\mathcal B}_{k}$ be the associated base subset of ${\mathcal G}_{k}$.
Denote by $s_{k}$ the cardinality of ${\mathcal B}_{k}$.

Clearly, we can suppose that $\sigma(i)=n+i$ for every $i\le n$
(see the definition of a symplectic base).
Then for each $i\in \{1,\dots,n\}$ every element of
${\mathcal B}_{n-1}$ contains precisely one of the points $p_{i}$ or $p_{n+i}$.
This implies that
$$s_{n-1}=2^{n}.$$
If $k<n-1$ then each $U\in {\mathcal B}_{k+1}$ contains $k+2$
distinct elements of ${\mathcal B}_{k}$
and every $S\in {\mathcal B}_{k}$ is contained in $2(n-k-1)$ distinct
subspaces belonging to ${\mathcal B}_{k+1}$.
Thus
$$s_{k}=s_{k+1}\frac{k+2}{2(n-k-1)}.$$
Step by step we get
$$s_{k}=s_{n-1}\frac{n}{2}\times\frac{n-1}{2\cdot 2}\times\dots\times\frac{k+2}{2(n-k-1)}$$
which gives the claim.
\end{proof}

\begin{prop}\label{pr2}
For any two $k$-dimensional totally isotropic subspaces
there is a base subset of ${\mathcal G}_{k}$ containing them.
\end{prop}

Proposition \ref{pr2} can be obtained by an immediate verification
or can be drawn from the fact that for any two flags there is an
apartment containing them.

\section{Result}
From this moment we suppose that
$V$ and $V'$ are $2n$-dimensional vector space over fields $F$ and $F'$
(respectively) and
$$\Omega: V\times V\to F,\;\;\;\Omega': V'\times V'\to F'$$
are non-degenerate symplectic forms.
Let $\Pi=(P,{\mathcal L})$ and $\Pi'=(P',{\mathcal L}')$
be the $(2n-1)$-dimensional  projective spaces  associated with $V$ and $V'$,
respectively.

An injection $f:P\to P'$ is called an {\it embedding} of $\Pi$ to $\Pi'$
if it maps lines to subsets of lines and for any line
$L'\in {\mathcal L}'$ there is at most one line $L\in {\mathcal L}$
such that $f(L)\subset L'$.
An embedding is said to be {\it strong} if it sends independent subsets
to independent subsets.
Every strong embedding of $\Pi$ to $\Pi$ is
induced by a semilinear injection of $V$ to $V'$ preserving the
linear independence \cite{FaureFrolicher,Faure,Havlicek}.

Our projective spaces have the same dimension, and
strong embeddings of $\Pi$ to $\Pi'$ (if they exist) map bases to bases.
An example given in \cite{HuangKreuzer} shows that strong
embeddings of $\Pi$ to $\Pi'$ can not be characterized as mappings
sending bases of $\Pi$ to bases of $\Pi$.
However, if $f:P\to P'$ transfers symplectic bases to symplectic bases
then $f$ is a strong embedding of $\Pi$ to $\Pi'$ and
for any two points $p,q\in P$
$$
p\perp q\;\Longrightarrow f(p)\perp f(q)\;\mbox{ and }\;
p\not\perp q\;\Longrightarrow f(p)\not\perp f(q),
$$
see \cite{Pankov4}.
Since a surjective embedding is a collineation,
every surjection of $P$ to $P'$ sending
symplectic bases to symplectic bases
is a collineation of $\Pi$ to $\Pi'$ preserving the orthogonal relation.

In what follows embeddings and collineations sending symplectic bases to
symplectic bases will be called {\it symplectic}.

Denote by ${\mathcal G}_{k}$ and ${\mathcal G}'_{k}$
the sets of $k$-dimensional totally isotrop\-i\-c subspaces
of $\Pi$ and $\Pi'$, respectively.

Let $f:P\to P'$ be a symplectic embedding of $\Pi$ to $\Pi'$.
For each $S\in {\mathcal G}_{k}$ the subspace spanned by $f(S)$
is an element of ${\mathcal G}'_{k}$. The mapping
$$(f)_{k}:{\mathcal G}_{k}\to {\mathcal G}'_{k}$$
$$S\to \overline{f(S)}$$
(we write $\overline{X}$ for the subspace spanned by $X$)
is an injection sending base subsets to base subsets.
If $f$ is a collineation then every $(f)_{k}$ is bijective.
Conversely, an easy verification shows that
if $(f)_{k}$ is bijective for certain $k$ then
$f$ is a collineation.

\begin{theorem}\label{theorem}
If a mapping of ${\mathcal G}_{k}$ to ${\mathcal G}'_{k}$
transfers base subsets to base subsets
then it is induced by a symplectic embedding of $\Pi$ to $\Pi'$.
\end{theorem}

\begin{cor}
Every surjection ${\mathcal G}_{k}$ to ${\mathcal G}'_{k}$
sending base subsets to base subsets is induced by a symplectic
collineation of $\Pi$ to $\Pi'$.
\end{cor}

For $k=n-1$ these results were established in \cite{Pankov3}.
For $n=2$ they can be drawn from well-known properties
of generalized quadrangles \cite{VM}.

Our proof of  Theorem \ref{theorem} is based on elementary properties of
so-called inexact subsets (Section 4).
If $k=n-1$ then all maximal inexact subsets are of the same type.
The case when $k<n-1$ is more complicated: there are two different types of
maximal inexact subsets.

Two elements of ${\mathcal G}_{k}$ are called {\it adjacent} if their intersection
belongs to ${\mathcal G}_{k-1}$.
We say that two elements of ${\mathcal G}_{k}$ are {\it ortho-adjacent} if they
are orthogonal and adjacent; this is possible only if $k<n-1$.
Using inexact subsets we characterize the adjacency and ortho-adjacency relations
in terms of base subsets.
This characterization shows that every mapping
of ${\mathcal G}_{k}$ to ${\mathcal G}'_{k}$ sending base subsets to base subsets
is adjacency and ortho-adjacency preserving (Section 6);
after that arguments in the spirit of \cite{Chow} give the claim (Section 7).

\section{Inexact subsets}
Let $n\ge 3$ and $B=\{p_{1},\dots,p_{2n}\}$ be a symplectic base of $\Pi$.
Denote by
${\mathcal B}$ the base subset of ${\mathcal G}_{k}$ associated with $B$.
By the definition,
${\mathcal B}$ consists of all $k$-dimensional subspaces
$$\overline{\{p_{i_{1}},\dots, p_{i_{k+1}}\}}$$
such that
$$\{i_{1},\dots,i_{k+1}\}\cap \{\sigma(i_{1}),\dots,\sigma(i_{k+1})\}=\emptyset.$$
If $k=m-1$ then every element of
${\mathcal B}$ contains precisely one of the points $p_{i}$ or $p_{\sigma(i)}$
for each $i$.

We write ${\mathcal B}(+i)$ and ${\mathcal B}(-i)$ for the sets of all
elements of ${\mathcal B}$ which  contain $p_{i}$ or
do not contain $p_{i}$, respectively.
For any $i_{1},\dots, i_{s}$ and $j_{1},\dots, j_{u}$
belonging to $\{1,\dots,2n\}$
we define
$${\mathcal B}(+i_{1},\dots, +i_{s},-j_{1},\dots, -j_{u}):=
{\mathcal B}(+i_{1})\cap\dots\cap {\mathcal B}(+i_{s})\cap
{\mathcal B}(-j_{1})\cap \dots\cap{\mathcal B}(-j_{u}).$$
The set of all elements of ${\mathcal B}$ incident with a subspace $S$
will be denoted by ${\mathcal B}(S)$ (this set may be empty).
Then ${\mathcal B}(-i)$ coincides with
${\mathcal B}(S)$, where $S$ is the subspace spanned by $B\setminus \{p_{i}\}$.
It is trivial that
$${\mathcal B}(+i)={\mathcal B}(+i,-\sigma(i))$$
and for the case when $k=m-1$ we have
$${\mathcal B}(-i)={\mathcal B}(+\sigma(i))={\mathcal B}(+\sigma(i),-i).$$

Let ${\mathcal R}\subset {\mathcal B}$.
We say that ${\mathcal R}$ is {\it exact} if
there is only one base subset of ${\mathcal G}_{k}$ containing ${\mathcal R}$;
otherwise, ${\mathcal R}$ will be called {\it inexact}.
If ${\mathcal R}\cap {\mathcal B}(+i)$ is not empty then we define
$S_{i}({\mathcal R})$ as the intersection of all subspaces
belonging to ${\mathcal R}$ and containing $p_{i}$,
and we define $S_{i}({\mathcal R}):=\emptyset$ if
the intersection of ${\mathcal R}$ and ${\mathcal B}(+i)$ is empty.
If
$$S_{i}({\mathcal R})=p_{i}$$
for all $i$ then ${\mathcal R}$ is exact;
the converse fails.

\begin{lemma}\label{lemma1}
Let ${\mathcal R}\subset {\mathcal B}$.
Suppose that there exist $i,j$ such that $j\ne i,\sigma(i)$
and
$$p_{j}\in S_{i}({\mathcal R}),\;\;
p_{\sigma(i)}\in S_{\sigma(j)}({\mathcal R}).$$
Then  ${\mathcal R}$ is inexact.
\end{lemma}

\begin{proof}
On the line $p_{i}p_{j}$ we choose a point $p'_{i}$
different from $p_{i}$ and $p_{j}$.
The line $p_{\sigma(i)}p_{\sigma(j)}$ contains a unique point orthogonal to $p'_{i}$;
we denote this point by $p'_{\sigma(j)}$.
Then
$$(B\setminus\{p_{i},p_{\sigma(j)}\})\cup \{p'_{i},p'_{\sigma(j)}\}$$
is a symplectic base.
The associated base subset of ${\mathcal G}_{k}$
contains ${\mathcal R}$ and we get the claim.
\end{proof}

\begin{prop}\label{pr3}
The subset ${\mathcal B}(-i)$ is inexact;
moreover, if $k<n-1$ then this is a maximal inexact subset.
For the case when $k=n-1$ the inexact subset ${\mathcal B}(-i)$ is not maximal.
\end{prop}

\begin{proof}
Let us  take a point $p'_{i}$ on the line $p_{i}p_{\sigma(i)}$
different from $p_{i}$ and $p_{\sigma(i)}$.
Then
$$(B\setminus\{p_{i}\})\cup \{p'_{i}\}$$
is a symplectic base and the associated base subset of
${\mathcal G}_{k}$ contains ${\mathcal B}(-i)$.
Hence this subset is inexact.

Let $k<n-1$.
For any $j\ne i$ we can choose distinct
$$i_{1},\dots,i_{k}\in \{1,\dots,2n\}\setminus\{i,j,\sigma(i),\sigma(j)\}$$
such that
$$\{i_{1},\dots,i_{k}\}\cap \{\sigma(i_{1}),\dots,\sigma(i_{k})\}=\emptyset.$$
The subspaces spanned by
$$p_{i_{1}},\dots,p_{i_{k}},p_{j}\;\mbox{ and }\;
p_{\sigma(i_{1})},\dots,p_{\sigma(i_{k})},p_{j}$$
belong to ${\mathcal B}(-i)$.
Since the intersection of these subspaces is $p_{j}$, we have
\begin{equation}\label{eq-1}
S_{j}({\mathcal B}(-i))=p_{j}\;\mbox{ if }\; j\ne i.
\end{equation}
Let $U$ be an arbitrary taken element of
$${\mathcal B}\setminus {\mathcal B}(-i)={\mathcal B}(+i).$$
This subspace is spanned by $p_{i}$ and some $p_{i_{1}},\dots,p_{i_{k}}$.
Since $p_{i}$ is a unique point of $U$ orthogonal to
$p_{\sigma(i_{1})},\dots,p_{\sigma(i_{k})}$,
\eqref{eq-1} shows that the subset
\begin{equation}\label{eq-2}
{\mathcal B}(-i)\cup \{U\}
\end{equation}
is exact.
This implies that the inexact subset ${\mathcal B}(-i)$ is maximal.

Now let $k=n-1$.
We take an arbitrary element $U\in {\mathcal B}(+i)$.
There exists $j$ such that $p_{\sigma(j)}$ does not belongs to $U$.
Then $p_{j}$ is a point of the subspace
$$S_{i}({\mathcal B}(-i)\cup \{U\})=U.$$
Since $p_{\sigma(i)}$ belongs to every element of ${\mathcal B}(-i)$ and
$p_{\sigma(j)}$ does not belongs to $U$,
$$
S_{\sigma(j)}({\mathcal B}(-i))=S_{\sigma(j)}({\mathcal B}(-i)\cup \{U\})
$$
contains $p_{\sigma(i)}$.
By Lemma \ref{lemma1}, the subset \eqref{eq-2} is inexact
and the inexact subset ${\mathcal B}(-i)$ is not maximal.
\end{proof}

\begin{prop}\label{pr4}
If $j\ne i,\sigma(i)$ then
$${\mathcal R}_{ij}:={\mathcal B}(+i,+j)\cup{\mathcal B}(+\sigma(i),+\sigma(j))
\cup{\mathcal B}(-i,-\sigma(j))$$
is a maximal inexact subset.
\end{prop}

If $k=n-1$ then
$${\mathcal R}_{ij}={\mathcal B}(+i,+j)\cup{\mathcal B}(-i).$$

\begin{proof}
Since
$$S_{i}({\mathcal R}_{ij})=p_{i}p_{j}\;\mbox{ and }\;
S_{\sigma(j)}({\mathcal R}_{ij})= p_{\sigma(j)}p_{\sigma(i)},$$
Lemma \ref{lemma1} shows that ${\mathcal R}_{ij}$ is inexact.
We want to show that
\begin{equation}\label{eq-3}
S_{l}({\mathcal R}_{ij})=p_{l}\;\mbox{ if }\; l\ne i,\sigma(j).
\end{equation}

Let $l\ne i,j,\sigma(i),\sigma(j)$.
If $k\ge 2$ then there exists
$$i_{1},\dots,i_{k-2}\in \{1,\dots,n\}\setminus\{i,j,\sigma(i),\sigma(j),l,\sigma(l)\}$$
such that
$$\{i_{1},\dots,i_{k}\}\cap \{\sigma(i_{1}),\dots,\sigma(i_{k})\}=\emptyset;$$
the subspaces spanned by
$$p_{i_{1}},\dots,p_{i_{k-2}},p_{l},p_{i},p_{j}\;\mbox{ and }
p_{\sigma(i_{1})},\dots,p_{\sigma(i_{k-2})},p_{l},
p_{\sigma(i)},p_{\sigma(j)}$$
are elements of ${\mathcal R}_{ij}$ intersecting in the point $p_{l}$.
If $k=1$ then the lines $p_{l}p_{\sigma(i)}$ and $p_{l}p_{j}$ are as required.

Now we choose distinct
$$i_{1},\dots,i_{k-1}\in \{1,\dots,n\}\setminus\{i,j,\sigma(i),\sigma(j)\}$$
such that
$$\{i_{1},\dots,i_{k-1}\}\cap \{\sigma(i_{1}),\dots,\sigma(i_{k-1})\}=\emptyset$$
and consider the subspace spanned by
$$p_{i_{1}},\dots,p_{i_{k-2}},p_{j},p_{\sigma(i)}.$$
This subspace intersects
the subspaces spanned by
$$p_{i_{1}},\dots,p_{i_{k-1}},p_{j},p_{i}\;\mbox{ and }\;
p_{i_{1}},\dots,p_{i_{k-1}},p_{\sigma(i)},p_{\sigma(j)}$$
precisely  in the points $p_{j}$ and  $p_{\sigma(i)}$, respectively.
Since all these subspaces are elements of ${\mathcal R}_{ij}$, we get
\eqref{eq-3} for $l=j,\sigma(i)$.

A direct verification shows that
$${\mathcal B}\setminus {\mathcal R}_{ij}=
{\mathcal B}(+i,-j)\cup{\mathcal B}(+\sigma(j),-\sigma(i)).$$
Thus for every $U\in {\mathcal B}\setminus {\mathcal R}_{ij}$
one of the following possibilities is realized:
\begin{enumerate}
\item[(1)] $U\in {\mathcal B}(+i,-j)$ intersects $S_{i}({\mathcal R}_{ij})=p_{i}p_{j}$
by $p_{i}$,
\item[(2)] $U\in {\mathcal B}(+\sigma(j),-\sigma(i))$
intersects $S_{\sigma(j)}({\mathcal R}_{ij})= p_{\sigma(j)}p_{\sigma(i)}$
by $p_{\sigma(j)}$.
\end{enumerate}
Since $p_{\sigma(j)}$ is a unique point of the line
$p_{\sigma(j)}p_{\sigma(i)}$ orthogonal to $p_{i}$
and $p_{i}$ is a unique point on $p_{i}p_{j}$ orthogonal to $p_{\sigma(j)}$,
the subset
$${\mathcal R}_{ij}\cup \{U\}$$
is exact for each $U$ belonging to ${\mathcal B}\setminus {\mathcal R}_{ij}$.
Thus the inexact subset ${\mathcal R}_{ij}$ is maximal.
\end{proof}

The maximal inexact subsets considered in Propositions \ref{pr3} and \ref{pr4}
will be called of the {\it first} and the {\it second} types, respectively.

\begin{prop}\label{pr5}
Every maximal inexact subset
is of the first or the second type.
In particular, if $k=n-1$ then each maximal inexact subset
is of the second type.
\end{prop}

\begin{proof}
Let ${\mathcal R}$ be a maximal inexect subset of ${\mathcal B}$,
and let $B'$ be another symplectic base of $\Pi$ such that
the associated base subset of ${\mathcal G}_{k}$ contains ${\mathcal R}$.
If certain $S_{i}({\mathcal R})$ is empty
then ${\mathcal R}\subset {\mathcal B}(-i)$.
If $k<n-1$ then the inverse inclusion holds
(since our inexact subset is maximal).
If $k=n-1$ then
$${\mathcal R}\subset {\mathcal B}(-i)\subset {\mathcal R}_{ij}$$
and ${\mathcal R}={\mathcal R}_{ij}$.

Now suppose that each $S_{i}({\mathcal R})$ is not empty.
Denote by $I$ the set of all $i$ such that
the dimension of $S_{i}({\mathcal R})$ is non-zero.
We take arbitrary $i\in I$ and suppose that
$S_{i}({\mathcal R})$ is spanned by $p_{i}$ and $p_{j_{1}},\dots, p_{j_{u}}$.
If the subspaces
$$M_{1}:=S_{\sigma(j_{1})}({\mathcal R}),\dots,
M_{u}:=S_{\sigma(j_{u})}({\mathcal R})$$
do not contain $p_{\sigma(i)}$ then
$p_{i}$ belongs to $M^{\perp}_{1},\dots,M^{\perp}_{u}$;
on the other hand
$$p_{j_{1}}\not\in M^{\perp}_{1},\dots,p_{j_{u}}\not\in M^{\perp}_{u}$$
and we have
$$M^{\perp}_{1}\cap\dots\cap M^{\perp}_{u}\cap S_{i}({\mathcal R})=p_{i};$$
since all $S_{l}({\mathcal R})$ and their orthogonal complements
are spanned by points of the base $B'$,
the point$p_{i}$ belongs to $B'$.
Therefore, there exist $i\in I$ and $j\ne i,\sigma(i)$
such that
$$p_{j}\in S_{i}({\mathcal R})\;\mbox{ and }\;p_{\sigma(i)}\in S_{\sigma(j)}({\mathcal R}).$$
Then ${\mathcal R}={\mathcal R}_{ij}$.
\end{proof}

Maximal inexact subsets of the same type have the same cardinality.
These cardinalities will be denoted by $c_{1}(k)$ and $c_{2}(k)$, respectively.
An immediate verification shows that each of the following possibilities
$$c_{1}(k)=c_{2}(k),\;c_{1}(k)<c_{2}(k),\; c_{1}(k)>c_{2}(k)$$
is realized.

\section{Complement subsets}
Let ${\mathcal B}$ be as in the previous section.
We say that ${\mathcal R}\subset {\mathcal B}$ is a {\it complement subset}
if ${\mathcal B}\setminus {\mathcal R}$
is a maximal inexact subset.
A complement subset is said to be of the {\it first} or the {\it second} type
if the corresponding maximal inexact subset is of the first or the second type,
respectively.
The compliment subsets for the maximal inexact subsets
from Propositions \ref{pr3} and \ref{pr4} are
$${\mathcal B}(+i)\;\mbox{ and }\;
{\mathcal B}(+i,-j)\cup{\mathcal B}(+\sigma(j),-\sigma(i)).$$
If $k=n-1$ then the second subset coincides with
$${\mathcal B}(+i,+\sigma(j))={\mathcal B}(+i,+\sigma(j),-j,-\sigma(i)).$$

\begin{lemma}\label{char-ad}
Let $k=n-1$.
Then $S,U\in {\mathcal B}$ are adjacent if and only if there are
precisely $\binom{k}{2}$ distinct complement subsets of ${\mathcal B}$ containing
both $S$ and $U$.
\end{lemma}

\begin{proof}
Denote by $m$ the dimension of $S\cap U$.
The complement subset ${\mathcal B}(+i,+j)$ contains our subspaces if and only if
$p_{i},p_{j}$ belong to $S\cap U$. Thus there are precisely $\binom{m+1}{2}$
distinct complement subsets of ${\mathcal B}$ containing $S$ and $U$.
\end{proof}

\begin{lemma}\label{max-in-4}
Let $k<n-1$ and ${\mathcal R}$ be a complement subset of ${\mathcal B}$.
If ${\mathcal R}$ is of the first type then there are precisely $4n-3$ distinct complement
subsets of ${\mathcal B}$ which do not intersect ${\mathcal R}$.
If ${\mathcal R}$ is
of the second type then there are precisely $4$ distinct complement
subsets of ${\mathcal B}$ which do not intersect ${\mathcal R}$.
\end{lemma}

To prove Lemma \ref{max-in-4} we use the following.

\begin{lemma}\label{spec-lemma}
Let $i,i',j,j'$ be elements of $\{1,\dots,2n\}$
such that $i\ne j$ and $i'\ne j'$.
If the intersection of
$${\mathcal B}(+i,-j)\;\mbox{ and }\;{\mathcal B}(+i',-j')$$
is empty then one of the following possibilities is realized:
$i'=\sigma(i)$, $i'=j$, $j'=i$.
\end{lemma}

\begin{proof}
Direct verification.
\end{proof}

\begin{proof}[Proof of Lemma \ref{max-in-4}]
Let us fix $l\in \{1,\dots,2n\}$ and consider
the complement subset ${\mathcal B}(+l)$.
If ${\mathcal B}(+i)$ is disjoint with ${\mathcal B}(+l)$
then $i=\sigma(l)$.
If for some $i,j\in \{1,\dots,2n\}$
the complement subset
$${\mathcal B}(+i,-j)\cup{\mathcal B}(+\sigma(j),-\sigma(i))$$
does not intersect ${\mathcal B}(+l)$ then one of the following possibilities is
realized:
\begin{enumerate}
\item[(1)] $i=\sigma(l)$, the condition $j\ne i,\sigma(i)$ shows that there are exactly
$2n-2$ possibilities for $j$;
\item[(2)] $j=l$ and there are exactly $2n-2$ possibilities for $i$
(since $i\ne j,\sigma(j)$).
\end{enumerate}
Now fix $i,j\in \{1,\dots,2n\}$ such that $j\ne i,\sigma(i)$
and consider the associated complement subset
\begin{equation}\label{compset}
{\mathcal B}(+i,-j)\cup{\mathcal B}(+\sigma(j),-\sigma(i)).
\end{equation}
There are only two complement subsets of the first type disjoint with
\eqref{compset}:
$${\mathcal B}(+\sigma(i))\;\mbox{ and }\;{\mathcal B}(+j).$$
If
$${\mathcal B}(+i',-j')\cup{\mathcal B}(+\sigma(j'),-\sigma(i'))$$
does not intersect \eqref{compset} then one of the following two possibilities
is realized:
$$i'=j,\;j'=i\;\mbox{ or }\;i'=\sigma(i), j'=\sigma(j)$$
(see Lemma \ref{spec-lemma}).
\end{proof}

\section{Main Lemma}
Let $f:{\mathcal G}_{k}\to {\mathcal G}'_{k}$
be a mapping which sends base subsets to base subsets.
Since for any two elements of ${\mathcal G}_{k}$
there exists a base subset containing them (Proposition \ref{pr2})
and the restriction of $f$ to every base subset of ${\mathcal G}_{k}$
is a bijection to a base subset of ${\mathcal G}'_{k}$,
the mapping $f$ is injective.

In this section the following statement will be proved.

\begin{lemma}[Main Lemma]\label{mainlemma}
Let $S,U\in {\mathcal G}_{k}$.
Then $S$ and $U$ are adjacent if and only if $f(S)$ and $f(U)$ are adjacent.
Moreover, for the case when $k<n-1$ the subspaces
$S$ and $U$ are ortho-adjacent if and only if
the same holds for $f(S)$ and $f(U)$.
\end{lemma}

Let ${\mathcal B}$ be a base subset of ${\mathcal G}_{k}$ containing $S$ and $U$.
Then ${\mathcal B}':=f({\mathcal B})$ is a base subset of ${\mathcal G}_{k}(\Omega')$
and the restriction $f|_{{\mathcal B}}$ is a bijection to ${\mathcal B}'$.

\begin{lemma}\label{max-in-2}
A subset ${\mathcal R}\subset {\mathcal B}$ is inexact
if and only if $f({\mathcal R})$ is inexact;
moreover,
${\mathcal R}$ is a maximal inexact subset if and only if
the same holds for $f({\mathcal R})$.
\end{lemma}

\begin{proof}
If ${\mathcal R}$ is inexact then there are two distinct
base subsets of ${\mathcal G}_{k}$ containing ${\mathcal R}$ and
their $f$-images are distinct base subsets of ${\mathcal G}'_{k}$
containing $f({\mathcal R})$, hence $f({\mathcal R})$ is inexact.
The base subsets ${\mathcal B}$ and ${\mathcal B}'$
have the same number of inexact subsets
and the first part of our statement is proved.

Now let ${\mathcal R}$ be a maximal inexact subset.

Suppose that $c_{1}(k)=c_{2}(k)$.
Then $f({\mathcal R})$ is an inexact subset of ${\mathcal B}'$
consisting of $c_{1}(k)=c_{2}(k)$ elements,
this inexact subset is maximal.
Since ${\mathcal B}$ and ${\mathcal B}'$
have the same number of maximal inexact subsets,
every maximal inexact subset of ${\mathcal B}'$
is the image of a maximal inexect subset of ${\mathcal B}$.

Consider the case when $c_{1}(k)>c_{2}(k)$
(the case $c_{1}(k)<c_{2}(k)$ is similar).
If ${\mathcal R}$ is of the first type
then the inexect subset $f({\mathcal R})$ consists of
$c_{1}(k)$ elements and the inequality $c_{1}(k)>c_{2}(k)$
guarantees that this is a maximal inexact subset of the first type.
The base subsets ${\mathcal B}$ and ${\mathcal B}'$
have the same number of maximal inexact subsets of the first type;
thus ${\mathcal R}$ is a maximal inexact subset of the first type
if and only if the same holds for $f({\mathcal R})$.
For the case when ${\mathcal R}$ is of the second type and
the inexact subset $f({\mathcal R})$ is not maximal,
we take a maximal inexact subset
${\mathcal R}'\subset {\mathcal B}'$
containing $f({\mathcal R})$; since
$$|{\mathcal R}'|>|f({\mathcal R})|=c_{2}(k),$$
${\mathcal R}'$ is of the first type and ${\mathcal R}$
is contained in the maximal inexact subset $f^{-1}({\mathcal R}')$;
the latter is impossible.
Similarly, we show that
every maximal inexact subset ${\mathcal R}'\subset {\mathcal B}'$
of the second type is the image of a maximal inexect subset of ${\mathcal B}$.
\end{proof}

\begin{lemma}\label{max-in-2,5}
${\mathcal R}\subset {\mathcal B}$ is a complement subset
if and only if $f({\mathcal R})$ is a complement subset of ${\mathcal B}'$.
\end{lemma}

\begin{proof}
This is a simple consequence of the previous lemma.
\end{proof}

For $k=n-1$ Main Lemma (Lemma \ref{mainlemma})
can be drawn directly from Lemmas \ref{char-ad} and \ref{max-in-2,5}.
In \cite{Pankov3} this statement was proved by a more complicated way.

\begin{lemma}\label{max-in-3}
The mapping $f|_{{\mathcal B}}$ together with the inverse mapping
preserve types of maximal inexact and complement subsets.
\end{lemma}

\begin{proof}
This statement is trivial if $k=n-1$ or $c_{1}(k)\ne c_{2}(k)$.
For the general case this follows from Lemma \ref{max-in-4}.
\end{proof}

We write ${\mathcal X}_{i}$ and ${\mathcal X}'_{i}$
for the sets of all $i$-dimensional subspaces
spanned by points of the symplectic bases
associated with ${\mathcal B}$ and ${\mathcal B}'$, respectively.

\begin{lemma}\label{max-in-5}
There exists a bijection $g:{\mathcal X}_{k+1}\to {\mathcal X}'_{k+1}$
such that
$$f({\mathcal B}(N))={\mathcal B}'(g(N))$$
for every $N\in {\mathcal X}_{k+1}$.
\end{lemma}

\begin{proof}
Lemma \ref{max-in-3} guarantees that
$f|_{{\mathcal B}}$ and the inverse mapping send
maximal inexact subsets of the first type to maximal inexact subsets of the first type.
This implies the existence of a bijection
$h:{\mathcal X}_{2n-2}\to {\mathcal X}'_{2n-2}$
such that
$$f({\mathcal B}(M))={\mathcal B}'(h(M))$$
for all $M\in {\mathcal X}_{2n-2}$.
Each $N\in {\mathcal X}_{k+1}$ can be presented
as the intersection of
$$M_{1},\dots, M_{2n-k-2}\in {\mathcal X}_{2n-2}.$$
Then
$$g(N):=\bigcap^{2n-k-2}_{i=1}\,h(M_{i})$$
is as required.
\end{proof}

Now we prove Lemma \ref{mainlemma} for $k<n-1$.
Two subspaces $S,U\in {\mathcal B}$ are adjacent
if and only if they belong to ${\mathcal B}(T)$ for certain $T\in {\mathcal X}_{k+1}$;
moreover, $S$ and $U$ are ortho-adjacent if and only if ${\mathcal B}(T)$ consists of
$k+2$ elements (in other words, $T$ is totally isotropic).
The required statement follows from Lemma \ref{max-in-5}.

\section{Proof of Theorem \ref{theorem}}

Let $M,N$ be a pair of incident subspaces of $\Pi$ such that
$\dim M<k<\dim N$.
We put $[M,N]_{k}$ for the set of $k$-dimensional subspaces of $\Pi$
incident with both $M$ and $N$;
for the case when $M=\emptyset$ or $N=P$ we write $(N]_{k}$ or $[M)_{k}$, respectively.

We say that ${\mathcal X}\subset {\mathcal G}_{k}$
is an $A$-{\it subset} if any two distinct elements of ${\mathcal X}$
are adjacent.

\begin{exmp}\label{ex1}{\rm
If $k<n-1$ and $N$ is an element of ${\mathcal G}_{k+1}$ then
$(N]_{k}$ is a maximal $A$-subset of ${\mathcal G}_{k}$.
Subsets of such type will be called {\it tops}.
Any two distinct elements of a top are ortho-adjacent.
}\end{exmp}

\begin{exmp}\label{ex2}{\rm
If $M$ belongs to ${\mathcal G}_{k-1}$ then
$$[M,M^{\perp}]_{k}=[M)_{k}\cap {\mathcal G}_{k}$$
is a maximal $A$-subset of ${\mathcal G}_{k}$.
Such maximal $A$-subsets are known as {\it stars},
they contain non-orthogonal elements.
}\end{exmp}

\begin{fact}[{\rm \cite{Chow},\cite{PPZ}}]
Each $A$-subset is contained in a maximal $A$-subset.
Every maximal $A$-subset of ${\mathcal G}_{n-1}$ is a star.
If $k<n-1$ then every maximal $A$-subset of ${\mathcal G}_{k}$
is a top or a star.
\end{fact}

The first part of Lemma \ref{mainlemma} says that $f$ transfers
$A$-subsets to $A$-subsets.
The second part of Lemma \ref{mainlemma} guarantees that
stars go to subsets of stars.
In other words, for any $M\in {\mathcal G}_{k-1}$
there exists $M'\in {\mathcal G}'_{k-1}$ such that
\begin{equation}\label{eq2.1}
f([M,M^{\perp}]_{k})\subset [M',M'^{\perp}]_{k}.
\end{equation}
Suppose that
$$f([M,M^{\perp}]_{k})\subset [M'',M''^{\perp}]_{k}$$
for other $M''\in {\mathcal G}'_{k-1}$.
Then $f([M,M^{\perp}]_{k})$ is contained in
the intersection of $[M',M'^{\perp}]_{k}$ and $[M'',M''^{\perp}]_{k}$.
This intersection is not empty only if
$M'=M''$ or $M'$ and $M''$ are ortho-adjacent;
but in the second case  our intersection consists of one element.
Thus there is  unique $M'\in {\mathcal G}'_{k-1}$ satisfying
\eqref{eq2.1}.
We have established the existence of a mapping
$$g:{\mathcal G}_{k-1}\to {\mathcal G}'_{k-1}$$
such that
$$f([M,M^{\perp}]_{k})\subset [g(M),g(M)^{\perp}]_{k}$$
for every $M\in {\mathcal G}_{k-1}$.
It is easy to see that
\begin{equation}\label{eq2.2}
g((N]_{k-1})\subset(f(N)]_{k-1}\;\;\;\;\;\forall\;N\in{\mathcal G}_{k}.
\end{equation}
Now we show that
{\it $g$ sends base subsets to base subsets}.

\begin{proof}
Let ${\mathcal B}_{k-1}$ be a base subset of ${\mathcal G}_{k-1}$
and $B$ be the associated symplectic base.
This  base defines a base subset ${\mathcal B}\subset {\mathcal G}_{k}$.
Now let $B'$ be the symplectic base associated with the base subset
${\mathcal B}':=f({\mathcal B})$
and ${\mathcal B}'_{k-1}$ be the base subset of ${\mathcal G}'_{k-1}$
defined by $B'$.
If $S\in {\mathcal B}_{k-1}$ then we take $U_{1},U_{2}\in {\mathcal B}$
such that $S=U_{1}\cap U_{2}$,
and \eqref{eq2.2} shows that
$$g(S)=f(U_{1})\cap f(U_{2})\in{\mathcal B}_{k-1}.$$
Thus
$g({\mathcal B}_{k-1})$ is contained in ${\mathcal B}'_{k-1}$.
Suppose that $g({\mathcal B}_{k-1})$ is a proper subset of ${\mathcal B}'_{k-1}$.
Then $g(S)=g(U)$ for some distinct $S,U\in {\mathcal B}_{k-1}$.
The $f$-image of
$${\mathcal B}(S)={\mathcal B}\cap[S,S^{\perp}]_{k}$$
is contained in
$${\mathcal B}'(g(S))={\mathcal B}'\cap [g(S),g(S)^{\perp}]_{k}.$$
Since these sets have the same cardinality,
$$f({\mathcal B}(S))={\mathcal B}'(g(S)).$$
Similarly,
$$f({\mathcal B}(U))={\mathcal B}'(g(U)).$$
The equality $f({\mathcal B}(S))=f({\mathcal B}(U))$
contradicts the injectivity of $f$.
Hence $g({\mathcal B}_{k-1})$ coincides with ${\mathcal B}'_{k-1}$.
\end{proof}

If $k=1$ then the mapping $g:P\to P'$ sends symplectic bases to symplectic bases,
by \cite{Pankov4} $g$ is a symplectic embedding of $\Pi$ to $\Pi'$,
and we have $f=(g)_{1}$.

Now suppose that $k>1$ and $g$ is induced by a symplectic embedding $h$
of $\Pi$ to $\Pi'$.
Let us consider an arbitrary element $S\in {\mathcal G}_{k}$ and
take ortho-adjacent $M,N\in{\mathcal G}_{k-1}$
such that $S=\overline{M\cup N}$.
Then
$$\{S\}=[M,M^{\perp}]_{k}\cap[N,N^{\perp}]_{k}$$
and $f(S)$ belongs to
the intersection of $[g(M),g(M)^{\perp}]_{k}$ and $[g(N),g(N)^{\perp}]_{k}$.
Since
$$g(M)=\overline{h(M)}\;\mbox{ and }\;g(N)=\overline{h(N)}$$
are ortho-adjacent,
the intersection of $[g(M),g(M)^{\perp}]_{k}$ and $[g(N),g(N)^{\perp}]_{k}$
consists of one element and
we have
$$f(S)=\overline{\overline{h(M)}\cup\overline{h(N)}}=\overline{h(S)}.$$
This means that $f$ is induced by $h$.
Therefore Theorem \ref{theorem} can be proved by induction.


\begin{thebibliography}{99}

\bibitem{AVM}
Abramenko P. and Van Maldeghem H.,
{\it On oppositions in spherical buildings and twin buildings},
Ann. Combinatorics 4(2000), 125-137.

\bibitem{Chow}
Chow W. L., {\it On the geometry of algebraic homogeneous spaces},
Ann. of Math. 50 (1949) 32-67.

\bibitem{Cohen}
Cohen A. M.,  {\it Point-Line Spaces Related to Buildings},
Handbook of incidence geometry.
Buildings and foundations (Edited by F. Buekenhout),
North-Holland, Amsterdam, 1995., p. 647-737.


\bibitem{FaureFrolicher} Faure C. A., Fr\"{o}licher A.,
{\it Morphisms of Projective Geometries and Semilinear maps},
Geom. Dedicata, 53 (1994), 237-262.

\bibitem{Faure} Faure C. A.,
{\it An Elementary Proof of the Fundamental Theorem of
Projective Geometry}, Geom. Dedicata, 90 (2002), 145-151.

\bibitem{Havlicek} Havlicek H.,
{\it A Generalization of Brauner's Theorem on Linear Mappings},
Mitt. Math. Sem. Univ. Gie$\beta$en, 215 (1994), 27-41.

\bibitem{HuangKreuzer}
Huang W.-l., Kreuzer A., {\it Basis preserving maps of linear spaces},
Arch. Math. (Basel) 64(1995), 6, 530-533.

\bibitem{Pankov1} Pankov M.,
Transformations of Grassmannians and automorphisms of classical groups,
J. Geom. 75 (2002), 132-150.

\bibitem{Pankov2} Pankov M.,
A characterization of geometrical mappings of Grassmann spaces,
Results Math. 45(2004), 319--327.

\bibitem{Pankov3} Pankov M.,
Mappings of the sets of invariant subspaces of null systems,
Beitr{\"a}ge zur Algebra Geom., 45(2004), 389-399.

\bibitem{Pankov4} Pankov M.,
Base preserving maps in symplectic geometry,
submitted to Abh. Math. Sem. Univ. Hamburg.

\bibitem{PPZ}
Pankov M., Pra\.{z}movski K., \.{Z}ynel M.,
{\it Geometry of polar Grassmann spaces}, accepted to Demonstratio Math.

\bibitem{Tits}
Tits J., {\it Buildings of spherical type and finite BN-pairs},
Lecture Notes in Mathematics 386, Springer, Berlin 1974.

\bibitem{VM} Van Maldeghem H., Private communication.
\end{thebibliography}
\end{document}